\documentclass[11pt]{article}
\usepackage{amsmath, amsthm}
\usepackage{latexsym}
\usepackage{graphicx}
\usepackage[scanall]{psfrag}
\makeatletter
\newfont{\footsc}{cmcsc10 at 8truept}
\newfont{\footbf}{cmbx10 at 8truept}
\newfont{\footrm}{cmr10 at 10truept}
\makeatother
\newtheorem{theorem}{\bf Theorem}
\newtheorem{proposition}{\bf Proposition}
\newtheorem{lemma}{\bf Lemma}
\newtheorem{corollary}{\bf Corollary}

\begin{document}
\title{Structural Properties of Bayesian Bandits with Exponential Family Distributions}

\author{Yaming Yu\\
\small Department of Statistics\\[-0.8ex]
\small University of California\\[-0.8ex] 
\small Irvine, CA 92697, USA\\[-0.8ex]
\small \texttt{yamingy@uci.edu}}

\date{}
\maketitle

\begin{abstract}
We study a bandit problem where observations from each arm have an exponential family distribution and different arms are assigned independent conjugate priors.  At each of $n$ stages, one arm is to be selected based on past observations.  The goal is to find a strategy that maximizes the expected discounted sum of the $n$ observations.  Two structural results hold in broad generality: (i) for a fixed prior weight, an arm becomes more desirable as its prior mean increases; (ii) for a fixed prior mean, an arm becomes more desirable as its prior weight decreases.  These  generalize and unify several results in the literature concerning specific problems including Bernoulli and normal bandits.  The second result captures an aspect of the exploration-exploitation dilemma in precise terms: given the same immediate payoff, the less one knows about an arm, the more desirable it becomes because there remains more information to be gained when selecting that arm.  For Bernoulli and normal bandits we also obtain extensions to nonconjugate priors. 


 
\end{abstract}

\section{Introduction}
At each of $n$ stages, an experimenter must take an observation from one of two stochastic processes (arms).  Let us adopt the Bayesian framework and assume that the experimenter's belief about an unknown arm is updated according to Bayes Theorem after each observation.  A strategy specifies which process to select at each stage.  The objective is to maximize the expected payoff, $\sum_{i=1}^n a_i Z_i$, where $Z_i$ is the observation at stage $i$ and $A_n\equiv (a_1, a_2,\ldots, a_n)$ is a discount sequence satisfying $a_i\geq 0$ and $\sum_{i=1}^n a_i >0$.  A strategy is optimal if it achieves the maximum expected payoff.  This is a finite-horizon two-armed bandit (Berry and Fristedt 1985), a classical problem in sequential decision theory. 

Bernoulli bandits, where each arm generates binary observations, are important as a model for clinical trials, and have received considerable attention (Berry 1972; Berry and Fristedt 1985).  Others such as normal (Chernoff 1968; Chernoff and Petkau 1986; Yao 2006) and Dirichlet bandits (Clayton and Berry 1985; Yu 2011) have also been extensively studied.  Bandit problems exhibit a well-known exploration-exploitation tradeoff.  Simply maximizing the immediate payoff is usually not an optimal strategy; one must allow for exploring an unknown arm for higher payoff later on.  From a Bayesian perspective, the optimal strategy is easily specified through backward induction, although its computation can be nontrivial.  If the discount sequence is geometric, then the problem reduces to several one-armed bandits (Gittins and Jones 1974; Gittins 1979; Whittle 1980; Kaspi and Mandelbaum 1998) and the optimal strategy is to choose an arm with the highest dynamic allocation index, or Gittins index.  Optimal strategies for general discount sequences are less tractable. 

The Gittins index possesses intriguing monotonicity properties with respect to prior specifications.  For example, Gittins and Wang (1992) show that the Gittins index decreases in $\tau>0$ for some special bandit arms: a Bernoulli arm whose unknown parameter has a ${\rm Beta}(\tau s, \tau(1-s))$ prior ($0<s<1$), or a normal arm whose unknown mean has a ${\rm N}(\mu, 1/\tau)$ prior ($\mu\in\mathbf{R}$).  In both cases $\tau$ is naturally interpreted as the amount of prior information.  Such monotonicity results therefore capture an aspect of the exploration-exploitation dilemma in precise terms: given the same immediate payoff, the less one knows about an arm, the more desirable it becomes since there is more room for exploration.  In the literature, however, this monotonicity is usually derived for one-armed bandits and on a case-by-case basis.  This paper aims to obtain more general results in a unified framework. 

The Bernoulli and normal bandits can be regarded as special cases of a general bandit where observations from each arm have an exponential family distribution.  Assume each arm is assigned an independent conjugate prior, which is characterized by a prior mean and a prior weight.  The prior mean specifies the immediate payoff of an arm, whereas the prior weight reflects the associated uncertainty.  For such problems we show that: (i) for fixed prior weight, the maximum expected payoff increases as the prior mean for any arm increases; (ii) for fixed prior mean, the maximum expected payoff increases as the prior weight for any arm decreases.  These generalize and unify several results in the literature concerning specific distributions.  Similar techniques yield parallel results for Dirichlet bandits, which do not fit in the one-parameter exponential family framework (Clayton and Berry 1985; Chattopadhyay 1994; Yu 2011). 

The rest of the paper is organized as follows.  After setting up the exponential family framework and introducing a few notions of stochastic ordering in Section~2, we present basic structural results such as a stay-on-a-winner rule in Section~3.  Section~4 contains the main results, including monotonicity of the value function with respect to prior weights.  Section~5 applies the results in Section~4 to one-armed bandits.  In particular, we show that the break-even value decreases as the prior weight of the unknown arm increases.  In Sections~6 and 7 we extend the monotonicity results to nonconjugate priors for Bernoulli and normal bandits, respectively.  Section~8 concludes with a brief discussion on an open problem. 

\section{Preliminaries}
Let $\nu$ be a $\sigma$-finite measure on $\mathbf{R}$ that is not a point mass.  Denote 
$$\psi(\theta)=\log \int e^{\theta x}\, {\rm d}\nu(x),\quad \theta\in \Theta,$$ 
where $\Theta$ is the natural parameter space defined as the set of $\theta\in \mathbf{R}$ such that $\psi(\theta)$ is finite.  We assume that $\Theta$ has a non-empty interior.  Suppose that given $\theta_i$, observations from arm $i$ are independent and identically distributed (i.i.d.) according to the density (relative to $\nu$) 
\begin{equation}
\label{density}
f(x|\theta_i)=e^{\theta_i x -\psi(\theta_i)}.
\end{equation}
Let us assume independent conjugate priors on $\theta_i,\ i=1,2,$ with Lebesgue density  \begin{equation}
\label{prior}
f(\theta_i|\gamma_i, \tau_i)\propto e^{\theta_i\gamma_i - \tau_i \psi(\theta_i)},\quad \theta_i\in \Theta.
\end{equation}
Let $\mathcal{K}$ denote the smallest open interval such that $\nu$ assigns no mass outside of the closure $\bar{\mathcal{K}}$.  To ensure that the priors are proper, we require $\tau_i>0$ and $\gamma_i/\tau_i\in \mathcal{K}$ (Brown 1986, Chapter 4).  As usual $\tau_i$ is regarded as the ``prior sample size'' and $\gamma_i$ the ``prior sum of observations''.
 We refer to (\ref{prior}) as the $(\gamma_i, \tau_i)$ prior and call this two-armed bandit with discount sequence $A_n$ the $(\gamma_1, \tau_1; \gamma_2, \tau_2; A_n)$ bandit.  Its value (i.e., maximum expected payoff) is denoted by $V(\gamma_1, \tau_1; \gamma_2, \tau_2; A_n)$. 

This framework unifies several well-studied bandit reward structures: (i) Bernoulli rewards whose unknown parameter has a ${\rm Beta}(\gamma, \tau -\gamma)$ prior; (ii) normal rewards whose unknown mean has a ${\rm N}(\gamma/\tau, 1/\tau)$ prior; (iii) exponential rewards whose unknown rate parameter has a ${\rm Gamma}(\tau+1, \gamma)$ prior; (iv) Poisson rewards whose unknown rate parameter has a ${\rm Gamma}(\gamma, \tau)$ prior.  Extensions to general priors for (i) and (ii) are considered in Sections 6 and 7, respectively. 

Let $V^i(\gamma_1,\tau_1; \gamma_2, \tau_2; A_n)$ be the expected payoff when selecting arm $i$ initially and using an optimal strategy thereafter.  Then 
\begin{equation}
\label{Vdef1}
V(\gamma_1, \tau_1;\gamma_2, \tau_2; A_n)=\max\left\{V^1(\gamma_1, \tau_1; \gamma_2, \tau_2; A_n), V^2(\gamma_1, \tau_1; \gamma_2, \tau_2; A_n)\right\}, 
\end{equation}
and it is optimal to start with the arm whose $V^i$ is larger.  Suppose arm $1$ is selected, resulting in an observation $X$.  By conjugacy, the posterior for $\theta_1$ is  again of the form of (\ref{prior}) with $(\gamma_1 + X, \tau_1 +1)$ in place of $(\gamma_1, \tau_1)$.  Thus we have 
\begin{align}
\label{Vdef2}
V^1(\gamma_1, \tau_1; \gamma_2, \tau_2; A_n) &= a_1 \mu_1 + \left. E\left[V(\gamma_1 + X, \tau_1 + 1; \gamma_2, \tau_2; A^1_n)\right|\gamma_1, \tau_1\right],\\
\label{Vdef3}
V^2(\gamma_1, \tau_1; \gamma_2, \tau_2; A_n) &= a_1 \mu_2 + \left. E\left[V(\gamma_1, \tau_1; \gamma_2 + Y, \tau_2 +1; A^1_n)\right|\gamma_2, \tau_2\right],
\end{align}
where $A_n^1=(a_2, a_3, \ldots, a_n)$ and $\mu_i$ denotes the expected value of an observation from arm $i$ under the $(\gamma_i, \tau_i)$ prior.  This $\mu_i$ is simply $\mu_i=\gamma_i/\tau_i$, which we refer to as the prior mean.  In $E[g(X)|\gamma_1, \tau_1]$, we use $X$ to denote a generic observation from arm 1 under the $(\gamma_1, \tau_1)$ prior; similarly for $Y$.  That is, the density of $X$ relative to $\nu$ is 
\begin{equation}
\label{margin}
f(x)\propto \int_\Theta e^{\theta (\gamma_1 + x) -(\tau_1+1) \psi(\theta)}\, {\rm d}\theta.
\end{equation}
The dynamic programming equations (\ref{Vdef1})--(\ref{Vdef3}) are crucial for both theoretical analysis and numerical computation of the optimal strategy. 

A key tool in our derivation is the notion of stochastic ordering (M\"{u}ller and Stoyan 2002; Shaked and Shanthikumar 2007).  We shall use the usual stochastic order $\leq_{\rm st}$, the convex order $\leq_{\rm cx}$, the likelihood ratio order $\leq_{\rm lr}$, and the relative log-concavity order $\leq_{\rm lc}$.  For random variables $Z_1$ and $Z_2$ taking values on $\mathbf{R}$, we write $Z_1\leq_{\rm st} Z_2$ (respectively, $Z_1\leq_{\rm cx} Z_2$), if $E\phi(Z_1)\leq E\phi(Z_2)$ for every increasing (respectively, convex) function $\phi$ such that the expectations exist.  If $Z_1\leq_{\rm st} Z_2$ then we also say $Z_2$ is to the right of $Z_1$.  If $Z_1$ and $Z_2$ have densities $f_1(z)$ and $f_2(z)$ respectively, supported on the same interval, then we write $Z_1\leq_{\rm lr} Z_2$ (respectively, $Z_1\leq_{\rm lc} Z_2$) if $\log \left(f_1(z)/f_2(z)\right)$ is decreasing (respectively, concave) in $z$.  For example, the $(\gamma, \tau)$ prior increases in the likelihood ratio order as $\gamma$ increases, and decreases in the relative log-concavity order as $\tau$ increases.  (We use $\leq_{\rm lr},\ \leq_{\rm st},\ \leq_{\rm lc}$ and $\leq_{\rm cx}$ with densities as well as random variables.)  Useful properties include the implication $\leq_{\rm lr} \Longrightarrow \leq_{\rm st}$.  Assuming equal means, it also holds that $\leq_{\rm lc}$ implies $\leq_{\rm cx}$.  Intuitively, the relative log-concavity order compares the amount of information as it is defined through curvatures of the log density functions.  Both $\leq_{\rm lr}$ and $\leq_{\rm lc}$ are preserved under the prior-to-posterior updating, which makes them ideal for studying structural properties in bandit problems.  The log-concavity order is also useful in other seemingly unrelated contexts (Whitt 1985; Yu 2009a, 2009b, 2010).  

\section{Stay-on-a-winner}
This section derives a basic monotonicity property of the optimal strategy: as the observation from an arm becomes larger, the inclination to pull that arm again also increases.  Under suitable conditions we prove a generalized stay-on-a-winner rule, which is a natural extension of the results for Bernoulli bandits (Bradt, Johnson and Karlin 1956; Berry 1972; Berry and Fristedt 1985). 

Let us define the advantage of arm 1 over arm 2 as
$$\Delta(\gamma_1, \tau_1; \gamma_2, \tau_2; A_n) =V^1(\gamma_1, \tau_1; \gamma_2, \tau_2; A_n) - V^2(\gamma_1, \tau_1; \gamma_2, \tau_2; A_n).$$ 
Define $\Delta^+=\max\{\Delta, 0\}$ and $\Delta^-=\min\{\Delta, 0\}$.  By considering the initial two pulls one can show (Berry 1972) 
\begin{align}
\label{diff1}
\Delta(\gamma_1, \tau_1; \gamma_2, \tau_2; A_n) = &(a_1-a_2)\left(\frac{\gamma_1}{\tau_1} -\frac{\gamma_2}{\tau_2}\right) \\ 
\label{diff2}
&+ E\left[\Delta^+(\gamma_1 +X, \tau_1+1; \gamma_2, \tau_2; A_n^1)|\gamma_1, \tau_1\right]\\
\label{diff3}
&+ E\left[\Delta^-(\gamma_1, \tau_1; \gamma_2 + Y, \tau_2+1; A_n^1)|\gamma_2, \tau_2\right].
\end{align}

Proposition~\ref{prop1} states that as the prior mean of arm 1 increases, so does the advantage of arm 1 over arm 2, assuming $A_n$ is decreasing.  This can be extended to non-conjugate priors.  Specifically, $\Delta$ increases as the prior for arm 1 becomes larger in the likelihood ratio order.  Extensions to general Markov decision problems are also possible (Rieder and Wagner 1991).  We provide a complete proof which serves as an introduction to the derivation of the main results in Section~4. 

\begin{proposition}
\label{prop1}
Suppose $A_n$ is decreasing.  Then $\Delta(\gamma_1, \tau_1; \gamma_2, \tau_2; A_n)$ increases in $\gamma_1$. 
\end{proposition}
\begin{proof}
The $n=1$ case is easy.  Let us use induction for $n\geq 2$.  In view of (\ref{diff1})--(\ref{diff3}), we only need to show that 
\begin{equation}
\label{diffmono1}
E\left[\Delta^+(\gamma_1 +X, \tau_1+1; \gamma_2, \tau_2; A_n^1)|\gamma_1, \tau_1\right]\quad {\rm and}
\end{equation}
\begin{equation}
\label{diffmono2}
E\left[\Delta^-(\gamma_1, \tau_1; \gamma_2 + Y, \tau_2+1; A_n^1)|\gamma_2, \tau_2\right]
\end{equation}
both increase in $\gamma_1$.  Monotonicity of (\ref{diffmono2}) follows from the induction hypothesis.  To handle (\ref{diffmono1}), let us consider $\gamma_1<\tilde{\gamma}_1$.  Let $\theta_1$ and $\tilde{\theta}_1$ have the $(\gamma_1, \tau_1)$ and $(\tilde{\gamma}_1, \tau_1)$ priors respectively.  Let $g(x)$ (respectively, $\tilde{g}(x)$) be the marginal density of $X$ if it is drawn according to (\ref{density}) given $\theta_1$ (respectively, $\tilde{\theta}_1$).  Note that $\theta_1\leq_{\rm lr} \tilde{\theta}_1$.  In view of (\ref{margin}), we know that $g\leq_{\rm lr} \tilde{g}$ by total positivity considerations (Karlin 1968, Chapter 3).  It follows that $g\leq_{\rm st} \tilde{g}$.  By the induction hypothesis, 
$$\phi(x)\equiv \Delta^+ (x, \tau_1 +1; \gamma_2, \tau_2; A_n^1)$$ 
increases in $x$.  Thus 
\begin{align}
\nonumber
E \left. \left[\phi(\gamma_1 + X) \right|\gamma_1, \tau_1 \right] & \leq \left. E\left[\phi(\tilde{\gamma}_1 + X) \right|\gamma_1, \tau_1 \right]\\
\label{diffmono4}
&\leq \left. E\left[\phi(\tilde{\gamma}_1 + X) \right|\tilde{\gamma}_1, \tau_1 \right],
\end{align}
where (\ref{diffmono4}) holds because $g\leq_{\rm st} \tilde{g}$.  Hence (\ref{diffmono1}) increases in $\gamma_1$.
\end{proof}

\begin{corollary}
Suppose $A_n$ is a decreasing sequence, and an observation $x$ is taken from arm 1 initially.  Then, at the second stage, either arm 1 is optimal for all $x$, or arm 2 is optimal for all $x$, or there exists some $x_*\in \mathcal{K}$ such that arm 1 is optimal if $x\geq x_*$ and arm 2 is optimal if $x\leq x_*$. 
\end{corollary}
\begin{proof}
We can show that $\Delta(\gamma_1+x, \tau_1+1; \gamma_2, \tau_2; A_n^1)$ is continuous in $x$.  (One method is to use the convexity result of Proposition~\ref{prop2} in Section 4.)  The claim then follows from Proposition~\ref{prop1}. 
\end{proof}

The next result, Theorem~\ref{thm0}, is a generalized stay-on-a-winner rule: under suitable conditions if an arm is optimal initially then it continues to be optimal at the next stage provided that the initial observation from that arm is large enough. 

\begin{theorem}
\label{thm0}
Assume $A_n$ is decreasing, $n\geq 2$, and either (i) $a_1=a_2$ or (ii) $\gamma_1/\tau_1\leq \gamma_2/\tau_2$ holds.  Assume $\Delta(\gamma_1, \tau_1; \gamma_2, \tau_2; A_n)\geq 0,$ i.e., arm 1 is optimal initially.  Then $\Delta(\gamma_1 + x, \tau_1 +1; \gamma_2, \tau_2; A_n^1)\geq 0$ for sufficiently large $x\in \bar{\mathcal{K}}$. 
\end{theorem}
\begin{proof}
We may assume $a_i>0$ for all $i\leq n$.  Let $U$ be the upper end point of $\mathcal{K}$.  If $U=\infty$, then using (\ref{diff1})--(\ref{diff3}), it is easy to show by induction that $\Delta(\gamma_1+x, \tau_1+1; \gamma_2, \tau_2; A_n^1)>0$ for sufficiently large $x$.  That is, the claim holds even without assuming that arm 1 is optimal initially.  Assume $U<\infty$ and $\Delta(\gamma_1, \tau_1; \gamma_2, \tau_2; A_n)\geq 0$.  By (\ref{diff1})--(\ref{diff3}) we have 
\begin{equation}
\label{bothzero}
0\leq E\left[\Delta^+(\gamma_1 +X, \tau_1+1; \gamma_2, \tau_2; A_n^1)|\gamma_1, \tau_1\right] + E\left[\Delta^-(\gamma_1, \tau_1; \gamma_2 + Y, \tau_2+1; A_n^1)|\gamma_2, \tau_2\right].
\end{equation}
Suppose the claim does not hold, i.e., $\Delta(\gamma_1 + x, \tau_1 +1; \gamma_2, \tau_2; A_n^1)< 0$ for all $x\in \mathcal{\bar{K}}$.  In particular, 
\begin{equation}
\label{del-}
\Delta(\gamma_1 + U, \tau_1 +1; \gamma_2, \tau_2; A_n^1)< 0.
\end{equation}  
Then it is necessary that both expectations in (\ref{bothzero}) are zero.  That is, 
$$\Delta(\gamma_1, \tau_1; \gamma_2 +y, \tau_2+1; A_n^1)\geq 0\quad {\rm for\ all\ } y\in \mathcal{K}.$$  
By continuity, $\Delta(\gamma_1, \tau_1; \gamma_2 +U, \tau_2+1; A_n^1)\geq 0$.  However, the $(\gamma_1+U, \tau_1+1)$ prior is larger than the $(\gamma_1, \tau_1)$ prior in the likelihood ratio order.  The argument of Proposition \ref{prop1} yields
\begin{align*}
\Delta(\gamma_1 + U, \tau_1 +1; \gamma_2, \tau_2; A_n^1) &\geq \Delta(\gamma_1, \tau_1; \gamma_2, \tau_2; A_n^1)\\
&\geq \Delta(\gamma_1, \tau_1; \gamma_2+U, \tau_2+1; A_n^1)\geq 0, 
\end{align*} 
which contradicts (\ref{del-}). 
\end{proof}

\section{Monotonicity}
Proposition \ref{prop2} shows that the maximum expected payoff is an increasing and convex function of the prior mean of any arm.  The convexity will be useful in proving Theorem \ref{thm1} concerning monotonicity with respect to the prior weight. 
\begin{proposition}
\label{prop2}
$V(\gamma_1, \tau_1; \gamma_2, \tau_2; A_n)$ is increasing and convex in each of $\gamma_i,\ i=1,2$.
\end{proposition}
\begin{proof}
Monotonicity holds by the same argument that proves Proposition \ref{prop1}.  Let us focus on the convexity with respect to $\gamma_1$.  The $n=1$ case is easy.  For $n\geq 2$ we use induction.  Note that by (\ref{Vdef1})--(\ref{Vdef3}) it suffices to show that both 
\begin{equation}
\label{convex1}
\left. E\left[V(\gamma_1 + X, \tau_1 +1; \gamma_2, \tau_2; A_n^1) \right|\gamma_1, \tau_1 \right]\quad {\rm and}
\end{equation}
\begin{equation}
\label{convex2}
\left. E\left[V(\gamma_1, \tau_1; \gamma_2 +Y, \tau_2 +1; A_n^1)\right|\gamma_2, \tau_2\right]
\end{equation}
are convex in $\gamma_1$.  The claim for (\ref{convex2}) follows from the induction hypothesis.  To deal with (\ref{convex1}), suppose $\gamma_1< \tilde{\gamma}_1$.  Denote the marginal of $X$ when the prior on $\theta$ is $(\gamma_1, \tau_1)$ (respectively, $(\tilde{\gamma}_1, \tau_1)$) by $g$ (respectively, $\tilde{g}$). 
Then $g\leq_{\rm st} \tilde{g}$ as in the proof of Proposition~\ref{prop1}.  By the induction hypothesis, 
$$\phi(x)\equiv V(x, \tau_1 +1; \gamma_2, \tau_2; A_n^1)$$ 
is convex in $x$.  Moreover, 
\begin{align}
\nonumber
E &\left. \left[\phi(\gamma_1 + X) \right|\gamma_1, \tau_1 \right] - E \left. \left[\phi\left(\frac{\gamma_1+\tilde{\gamma}_1}{2}+ X\right) \right|\gamma_1, \tau_1 \right]\\
\label{ineq1}
&\geq  
 E \left. \left[\eta(X)\right|\gamma_1, \tau_1 \right]\\
 \label{ineq2}
 &\geq  E \left. \left[\eta(X)\right|\tilde{\gamma}_1, \tau_1 \right]
\end{align}
where 
$$\eta(x) \equiv \phi\left(\frac{\gamma_1+\tilde{\gamma}_1}{2} + x\right) -\phi(\tilde{\gamma}_1 + x).$$
The inequality (\ref{ineq1}) holds because $\phi$ is convex; (\ref{ineq2}) holds because $\eta$ is decreasing and  $g\leq_{\rm st} \tilde{g}$.  Rearranging we get 
$$E\left. \left[\phi(\gamma_1 + X) \right|\gamma_1, \tau_1 \right] + E \left. \left[\phi(\tilde{\gamma}_1 + X) \right|\tilde{\gamma}_1, \tau_1 \right] \geq 2 E \phi\left(\frac{\gamma_1+\tilde{\gamma}_1}{2} + X^*\right)$$
where $X^*$ has the following distribution.  Given $\theta$, $X^*$ is distributed according to (\ref{density}); the prior on $\theta$ is a half-half mixture of $(\gamma_1, \tau_1)$ and $(\tilde{\gamma}_1, \tau_1)$.  Denote this mixture density by $h^*(\theta)$, and the $((\gamma_1 + \tilde{\gamma}_1)/2, \tau_1)$ prior density by $h(\theta)$.  Then $h(\theta)\leq_{\rm lc} h^*(\theta)$, because log-convexity is closed under mixtures (Marshall and Olkin 1979).  Consider the difference between the marginal densities
$$D(x)\equiv \int_\Theta e^{x\theta -\psi(\theta)} \left[h(\theta) - h^*(\theta)\right]\, {\rm d}\theta.$$
Relative log-concavity implies that, as $\theta$ traverses $\Theta,\ h(\theta)- h^*(\theta)$ changes signs at most twice and, in the case of two changes, the sign sequence is $-, +, -$.  By the variation-diminishing properties of the Laplace transform (Karlin 1968, Chapter 5), $D(x)$ has at most two changes of sign, and in the case of two changes, the sign sequence is $-, +, -$.  Note that, when the prior is either $h$ or $h^*$, the marginal mean of $X$ is the same, namely $(\gamma_1 +\tilde{\gamma}_1)/(2\tau_1)$.  Hence it is not possible for $D(x)$ to change signs exactly once.  Unless $D(x)\equiv 0$, its sign sequence must be $-, +, -$.  It follows that the marginal distribution of $X$ becomes larger in the convex order when $\tilde{h}(\theta)$ replaces $h(\theta)$ as the prior for $\theta$ (see, e.g., Yu 2010, Lemma 1).  Using the convexity of $\phi$ again, we obtain
$$E \phi\left(\frac{\gamma_1+\tilde{\gamma}_1}{2} + X^*\right)\geq 
\left. E \left[\phi\left(\frac{\gamma_1+\tilde{\gamma}_1}{2} + X\right)\right| \frac{\gamma_1+\tilde{\gamma}_1}{2}, \tau_1\right].
$$
It follows that $E\left[\phi(\gamma_1 +X)|\gamma_1, \tau_1\right]$, i.e., (\ref{convex1}), is convex in $\gamma_1$, as required. 
\end{proof}

Our main result, Theorem \ref{thm1}, shows that the value of the bandit decreases as the prior weight of  an arm increases.  That is, given the same immediate payoff, an arm becomes less desirable as the amount of information about it increases. 

\begin{theorem}
\label{thm1}
$V(c\gamma_1, c\tau_1; \gamma_2, \tau_2; A_n)$ decreases in $c\in (0, \infty)$. 
\end{theorem}
\begin{proof}
Let us use induction on $n$.  The $n=1$ case is easy.  Suppose $n\geq 2$.  In view of (\ref{Vdef1})--(\ref{Vdef3}), we only need to show that 
\begin{equation}
\label{decrease1}
\left. E\left[V(c\gamma_1 + X, c\tau_1 +1; \gamma_2, \tau_2; A_n^1) \right|c\gamma_1, c\tau_1 \right]\quad {\rm and}
\end{equation}
\begin{equation}
\label{decrease2}
\left. E\left[V(c\gamma_1, c\tau_1; \gamma_2 +Y, \tau_2 +1; A_n^1)\right|\gamma_2, \tau_2\right]
\end{equation}
both decrease in $c$.  By the induction hypothesis, (\ref{decrease2}) decreases in $c$.  To deal with (\ref{decrease1}), suppose $0< c<\tilde{c}$ and denote $\xi = (c\tau_1 + 1)/(\tilde{c} \tau_1 +1)$.  We get 
\begin{align}
\nonumber
& \left. E\left[V(c\gamma_1 + X, c\tau_1 +1; \gamma_2, \tau_2; A_n^1) \right|c\gamma_1, c\tau_1 \right]\\
\label{ineqA}
&\geq \left. E\left[V(\xi (\tilde{c}\gamma_1 + X), c\tau_1 +1; \gamma_2, \tau_2; A_n^1) \right|c\gamma_1, c\tau_1 \right]\\
\label{ineqB}
&\geq \left. E\left[V(\tilde{c}\gamma_1 + X, \tilde{c}\tau_1 +1; \gamma_2, \tau_2; A_n^1) \right|c\gamma_1, c\tau_1 \right]\\
\label{ineqC}
&\geq \left. E\left[V(\tilde{c}\gamma_1 + X, \tilde{c}\tau_1 +1; \gamma_2, \tau_2; A_n^1) \right|\tilde{c}\gamma_1, \tilde{c}\tau_1 \right].
\end{align}
The inequality (\ref{ineqA}) holds by the convexity of $V$ as shown by Proposition~\ref{prop2}, noting 
$$\xi (\tilde{c}\gamma_1 +X)\leq_{\rm cx} c\gamma_1 + X$$
(see Lemma~\ref{lem3} in Section~7, or Shaked and Shanthikumar 2007, Theorem 3.A.18).  The inequality (\ref{ineqB}) holds by the induction hypothesis, as $\xi<1$.  The inequality (\ref{ineqC}) holds by an argument similar to the proof of Proposition \ref{prop2}.  Specifically, the prior $(\tilde{c}\gamma_1, \tilde{c}\tau_1)$ is log-concave relative to $(c \gamma_1, c \tau_1)$.  Thus the marginal of $X$ increases in the convex order if $(c\gamma_1, c\tau_1)$ replaces $(\tilde{c}\gamma_1, \tilde{c}\tau_1)$ as the prior on $\theta$ (the mean of $X$ remains constant).  Overall  (\ref{decrease1}) decreases in $c$, as required. 
\end{proof}

{\bf Remark.}  Proposition~\ref{prop2} and Theorem~\ref{thm1} extend naturally to bandits with more than two arms.  We present the two-armed version for simplicity.  The discount sequence $A_n$ is only required to be nonnegative.  By approximation, this can be further extended to the infinite-horizon case assuming $\sum_{i=1}^\infty a_i <\infty$. 

\section{The one-armed case} 
This section considers the one-armed case assuming that arm 2 yields a constant payoff $\lambda$ at each pull.  We shall abuse the notation by calling this a $(\gamma, \tau; \lambda; A_n)$ bandit, where we drop the subscripts on $\gamma_1$ and $\tau_1$ for convenience.  Results in Section~4 are applied to derive monotonicity properties of the break-even value in this case.  It is also shown (Proposition~\ref{prop3}) that if both arms are optimal initially, then an observation from arm 1 that is less than its prior mean would make arm 2 optimal thereafter. 

A discount sequence $A_n=(a_1, a_2, \ldots)$ is called {\it regular} if, letting $b_j=\sum_{i\geq j} a_i$, we have 
$b_{j+1}^2\geq b_j b_{j+2}$ for all $j\geq 1$ (Berry and Fristedt 1979).  For regular discount sequences, our one-armed bandit is an optimal stopping problem, i.e., if at any stage the known arm becomes optimal then it remains optimal in all subsequent stages.  Moreover, if $A_n$ is regular and $a_1>0$, then there exists a break-even value $\Lambda(\gamma, \tau; A_n)$ for the $(\gamma, \tau; \lambda; A_n)$ bandit, such that arm 1 is optimal initially if and only if $\lambda \leq \Lambda(\gamma, \tau; A_n)$ and arm 2 is optimal initially if and only if $\lambda \geq \Lambda(\gamma, \tau; A_n)$.  For infinite-horizon geometric discounting, this break-even value is also known as the dynamic allocation index or Gittins index (Gittins and Jones 1974).  The following result holds by the optimal stopping characterization. 

\begin{lemma}
\label{lemlam}
If $A_n$ is regular and $a_1>0$, then $\Lambda(\gamma, \tau; A_n)$ is the smallest $\lambda$ such that 
$$V(\gamma, \tau; \lambda; A_n)\leq \lambda \sum_{i= 1}^n a_i .$$ 
\end{lemma}

Corollary~\ref{coro1} summarizes some monotonicity properties of $\Lambda(\gamma, \tau; A_n)$.  It extends to infinite-horizon regular discounting.  As special cases we recover the results of Gittins and Wang (1992) on Bernoulli and normal bandits with geometric discounting; see also Yao (2006). 

\begin{corollary}
\label{coro1}
If $A_n$ is regular and $a_1>0$, then $\Lambda(c\gamma, c\tau; A_n)$ decreases in $c>0$ and strictly increases in $\gamma$. 
\end{corollary}
\begin{proof}
Monotonicity in $c$ follows from Theorem \ref{thm1} and Lemma \ref{lemlam}.  Monotonicity in $\gamma$ follows from Proposition \ref{prop2} and Lemma \ref{lemlam}.  To show strict monotonicity, let us set $c=1$ and assume that $\gamma, \tilde{\gamma}$ satisfy $\gamma< \tilde{\gamma}$ and 
$$\Lambda(\gamma, \tau; A_n)=\Lambda(\tilde{\gamma}, \tau; A_n)\equiv \lambda_*.$$ 
Then, as in the proof of Proposition \ref{prop1}, we get 
\begin{align*}
\lambda_* \sum_{i=1}^n a_i &= a_1 \frac{\gamma}{\tau} + \left. E\left[V(\gamma+X, \tau + 1; \lambda_*; A_n^1)\right|\gamma, \tau\right]\\
& < a_1\frac{\tilde{\gamma}}{\tau} +\left. E\left[V(\gamma +X, \tau + 1; \lambda_*; A_n^1)\right|\gamma, \tau\right]\\
& \leq a_1\frac{\tilde{\gamma}}{\tau} +\left. E\left[V(\tilde{\gamma}+X, \tau + 1; \lambda_*; A_n^1)\right|\tilde{\gamma}, \tau\right]\\
& = \lambda_* \sum_{i=1}^n a_i, 
\end{align*}
which is a contradiction. 
\end{proof}

For a regular and positive discount sequence $A_n$, Proposition \ref{prop3} shows that there exists a break-even observation $b(\gamma, \tau; A_n)$ for the $(\gamma, \tau; \lambda; A_n)$ bandit such that if both arms are optimal initially, and an observation $x$ is taken from arm 1, then arm 1 remains optimal if $x\geq b(\gamma, \tau; A_n)$ and arm 2 becomes optimal if $x\leq b(\gamma, \tau; A_n)$.  Moreover, this break-even observation is no smaller than $\gamma/\tau$, the prior mean. 

\begin{proposition} 
\label{prop3}
Suppose $A_n$ is regular, $n\geq 2,$ and $a_1, a_2>0$.  Then there exists a unique $b(\gamma, \tau; A_n)\in \mathcal{K}$ such that 
$b(\gamma, \tau; A_n)\geq \gamma/\tau$ and  
\begin{align}
\label{lambreak1}
\Lambda(\gamma, \tau; A_n) &\geq \Lambda(\gamma + x, \tau + 1; A_n^1),\quad {\rm if}\ x\leq b(\gamma, \tau; A_n);\\
\label{lambreak2}
\Lambda(\gamma, \tau; A_n) &\leq \Lambda(\gamma + x, \tau + 1; A_n^1),\quad {\rm if}\ x\geq b(\gamma, \tau; A_n).
\end{align}
\end{proposition}

To prove Proposition~\ref{prop3} we need a continuity lemma.  Its proof, taken from Clayton and Berry (1985), is included for completeness. 
\begin{lemma}
\label{cont}
Suppose $A_n$ is regular and $a_1>0$.  Then $\Lambda(\gamma, \tau; A_n)$ is continuous in $\gamma$. 
\end{lemma} 
\begin{proof}
Fix $\gamma_0$ and note that $\lambda=\Lambda(\gamma, \tau; A_n)$ is the unique root of 
$$V^1(\gamma, \tau; \lambda; A_n) - V^2(\gamma, \tau; \lambda; A_n)=0.$$
By continuity of $V^1$ and $V^2$, we have
\begin{align*}
0&=\lim_{\gamma\uparrow\gamma_0} \left[V^1(\gamma, \tau; \Lambda(\gamma, \tau; A_n); A_n) - V^2(\gamma, \tau; \Lambda(\gamma, \tau; A_n); A_n)\right]\\
&= V^1(\gamma_0, \tau; \lim_{\gamma\uparrow \gamma_0} \Lambda(\gamma, \tau; A_n); A_n) - V^2(\gamma_0, \tau; \lim_{\gamma\uparrow\gamma_0} \Lambda(\gamma, \tau; A_n); A_n). 
\end{align*}
By uniqueness of $\Lambda$, we have $\lim_{\gamma\uparrow\gamma_0} \Lambda(\gamma, \tau; A_n)=\Lambda(\gamma_0, \tau; A_n)$.  Similarly, the limit holds when $\gamma\downarrow\gamma_0$. 
\end{proof}

\begin{proof}[Proof of Proposition~\ref{prop3}]
Let $U$ be the upper end point of $\mathcal{K}$.  If $U=\infty$ then $\Lambda(\gamma +x, \tau+1; A_n^1)\to \infty$ as $x\to \infty$ (the expected payoff by always selecting arm 1 becomes arbitrarily large).  If $U<\infty$ then we can show $\Lambda(\gamma + U, \tau+1; A_n^1) > \Lambda(\gamma, \tau; A_n)$ as follows.  Assume the contrary and consider the $(\gamma, \tau; \lambda_*; A_n)$ bandit with $\lambda_*=\Lambda(\gamma+U, \tau+1; A_n^1)$.  We have  
$$\lambda_* \sum_{i=1}^n a_i \leq a_1\frac{\gamma}{\tau} + \left. E\left[V(\gamma+X, \tau +1; \lambda_*; A_n^1)\right|\gamma, \tau\right].$$
Since $\gamma/\tau\in \mathcal{K}$ and $\mathcal{K}$ is open, we have $\lambda_* \geq (\gamma + U)/(\tau+1)> \gamma/\tau$.  Thus   
\begin{align*}
 \lambda_* \sum_{i=2}^n a_i & < \left. E\left[V(\gamma+X, \tau +1; \lambda_*; A_n^1)\right|\gamma, \tau\right]\\
&\leq V(\gamma + U, \tau+1; \lambda_*; A_n^1)\\
& =\lambda_* \sum_{i=2}^n a_i, 
\end{align*}
which is a contradiction.  We also have
\begin{align*}
\Lambda(\gamma, \tau; A_n) &\geq \Lambda(\gamma, \tau; A_n^1)\\
&\geq \Lambda(\gamma+\gamma/\tau, \tau+1; A_n^1)
\end{align*}
where the first inequality holds by the optimal stopping characterization, and the second by Corollary \ref{coro1}.

By Lemma \ref{cont} and Corollary \ref{coro1}, $\Lambda(\gamma + x, \tau + 1; A_n^1)$ is continuous and strictly increasing in $x$.  By the mean value theorem, there exists a unique $b(\gamma, \tau; A_n)\in [\gamma/\tau, U)$ such that (\ref{lambreak1}) and (\ref{lambreak2}) hold. 
\end{proof}

It is tempting to conjecture that $b(\gamma, \tau; A_n)\geq \Lambda(\gamma, \tau; A_n)$, which gives a tighter bound since $\Lambda(\gamma, \tau; A_n)\geq \gamma/\tau$.  However, our methods are not yet strong enough to resolve this conjecture.  Clayton and Berry (1985) conjectured and Yu (2011) proved an analogous bound for Dirichlet bandits. 

\section{Bernoulli bandits with general priors}
As noted earlier, results based on likelihood ratio orders, such as those in Section~3, may extend to nonconjugate priors.  This section shows that Theorem~\ref{thm1} can also be extended this way, at least in the Bernoulli case. 

Given $p_i,\ i=1,2,$ let us assume that observations from arm $i$ are i.i.d.\ ${\rm Bernoulli}(p_i)$.  Priors on $p_i$ are independent with densities $f_i$ with respect to a $\sigma$-finite measure $G$ on $[0,1]$.  We shall denote the value of this Bernoulli bandit with discount sequence $A_n$ by $V_{\rm B}(f_1; f_2; A_n)$.  Let $\mu(f)$ denote the mean of any prior $f$, i.e., $\mu(f)=\int_{[0,1]} p f(p)\, {\rm d}G(p)$. 

\begin{theorem}
\label{thm2}
If $f_1\leq_{\rm lc} \tilde{f}_1$ and $\mu(f_1)= \mu(\tilde{f}_1)$, then $V_{\rm B}(f_1; f_2; A_n)\leq V_{\rm B}(\tilde{f}_1; f_2; A_n)$.
\end{theorem}

Note that the ${\rm Beta}(c\alpha, c\beta)$ prior ($c, \alpha, \beta>0$) decreases in the relative log-concavity order as $c$ increases.  Theorem~\ref{thm2} therefore recovers the Bernoulli case of Theorem~\ref{thm1} for conjugate priors.

Let $\Lambda_{\rm B}(f; A_n)$ denote the break-even value of a one armed Bernoulli bandit whose unknown arm has prior $f$.  We obtain Corollary \ref{coro3} as a consequence of Theorem~\ref{thm2} and Lemma~\ref{lemlam}.

\begin{corollary}
\label{coro3}
Assume $A_n$ is regular and $a_1>0$.  If $f\leq_{\rm lc} \tilde{f}$ and $\mu(f)= \mu(\tilde{f})$, then $\Lambda_{\rm B}(f; A_n)\leq \Lambda_{\rm B}(\tilde{f}; A_n)$. 
\end{corollary}

Herschkorn (1997) posed the problem of identifying a variability ordering between priors so that both $V_{\rm B}$ and $\Lambda_{\rm B}$ are monotonic with respect to it.  Theorem~\ref{thm2} and Corollary \ref{coro3} show that there is indeed such an ordering, namely the relative log-concavity order (assuming equal means).  A conjecture of Herschkorn (1997) states that Corollary~\ref{coro3} holds under the weaker assumption $f\leq_{\rm cx} \tilde{f}$.  This conjecture remains open.

\begin{proof}[Proof of Theorem~\ref{thm2}]
The $n=1$ case is easy.  For $n\geq 2$ we use induction.  The equations (\ref{Vdef1})--(\ref{Vdef3}) become 
\begin{align}
\nonumber
V_{\rm B}(f_1; f_2; A_n) &= \max\{ V_{\rm B}^1(f_1; f_2; A_n),\, V_{\rm B}^2(f_1; f_2; A_n)\};\\
\label{Vb1}
V_{\rm B}^1(f_1; f_2; A_n) &=\mu(f_1) (a_1 + V_{\rm B}(\sigma f_1; f_2; A_n^1)) + (1-\mu(f_1)) V_{\rm B}(\phi f_1; f_2; A_n^1);\\
\nonumber 
V_{\rm B}^2(f_1; f_2; A_n) & =\mu(f_2) (a_1 + V_{\rm B}(f_1; \sigma f_2; A_n^1)) + (1-\mu(f_2)) V_{\rm B}(f_1; \phi f_2; A_n^1).
\end{align}
We use $\sigma f$ (respectively, $\phi f$) to denote the posterior density after observing one success (respectively, one failure).  That is, 
$$(\sigma f)(p)= \frac{f(p) p}{\mu(f)};\quad (\phi f)(p)= \frac{f(p) (1-p)}{1-\mu(f)}.$$ 

Let us assume $\tilde{f}_1$ is nondegenerate.  Because $f_1\leq_{\rm lc} \tilde{f}_1$ and $\mu(f_1)=\mu(\tilde{f}_1)$ we have $f_1\leq_{\rm cx} \tilde{f}_1$ (see, e.g., Yu 2010, Theorem 12).  Thus 
$$\mu(f_1) \mu(\sigma f_1)=\int_{[0,1]} p^2 f_1(p)\, {\rm d}G(p) \leq \int_{[0,1]} p^2 \tilde{f}_1(p)\, {\rm d}G(p) = \mu(\sigma \tilde{f}_1)\mu(\tilde{f}_1),$$
yielding $\mu(\sigma f_1) \leq \mu(\sigma \tilde{f}_1)$.  Similarly, $\mu(\phi f_1) \geq \mu(\phi \tilde{f}_1).$  Define 
$$\epsilon^*=\frac{\mu(\sigma \tilde{f}_1) -\mu(\sigma f_1)}{\mu(\sigma \tilde{f}_1) -\mu(\phi \tilde{f}_1)};\quad \epsilon_*=\frac{\mu(\phi f_1) -\mu(\phi \tilde{f}_1)}{\mu(\sigma \tilde{f}_1) -\mu(\phi \tilde{f}_1)}.$$
Then $\epsilon^*, \epsilon_*\in [0, 1)$.  Define 
$$g^*=(1-\epsilon^*)\sigma \tilde{f}_1 +\epsilon^*\phi\tilde{f}_1;\quad g_*=\epsilon_* \sigma \tilde{f}_1 +(1-\epsilon_*)\phi\tilde{f}_1.$$
Convexity of $V_{\rm B}$ with respect to mixtures gives 
\begin{align*}
V_{\rm B}(g^*; f_2; A_n^1) &\leq (1-\epsilon^*) V_{\rm B}(\sigma \tilde{f}_1; f_2; A_n^1) +\epsilon^* V_{\rm B}(\phi\tilde{f}_1; f_2; A_n^1);\\
V_{\rm B}(g_*; f_2; A_n^1) &\leq \epsilon_* V_{\rm B}(\sigma \tilde{f}_1; f_2; A_n^1) +(1-\epsilon_*) V_{\rm B}(\phi\tilde{f}_1; f_2; A_n^1).
\end{align*}
Noting $\mu(f_1) \epsilon^*=(1-\mu(f_1))\epsilon_*$, we add $\mu(f_1)$ times the first inequality to $1-\mu(f_1)$ times the second and get
\begin{align}
\nonumber
\mu(f_1) V_{\rm B}(g^*; f_2; A_n^1) &+(1-\mu(f_1)) V_{\rm B}(g_*; f_2; A_n^1)\\
\label{Vb2}
\leq & \mu(f_1) V_{\rm B}(\sigma \tilde{f}_1; f_2; A_n^1) + (1-\mu(f_1)) V_{\rm B}(\phi\tilde{f}_1; f_2; A_n^1).
\end{align}
The density $g^*$ is simply
$$g^*(p)=\left[\frac{p(1-\epsilon^*)}{\mu(f_1)} + \frac{(1-p)\epsilon^*}{1-\mu(f_1)}\right] \tilde{f}_1(p).$$
It is easy to check $(1-\epsilon^*)/\mu(f_1)\geq \epsilon^*/(1-\mu(f_1))$, which leads to
$$\sigma f_1\leq_{\rm lc} \sigma \tilde{f}_1 \leq_{\rm lc} g^*.$$
Moreover, $\sigma f_1$ and $g^*$ have the same mean.  By the induction hypothesis, we have
\begin{align}
\label{Vb3}
V_{\rm B}(\sigma f_1; f_2; A_n^1) &\leq V_{\rm B}(g^*; f_2; A_n^1).
\end{align}
Similarly, 
\begin{align}
\label{Vb4}
V_{\rm B}(\phi f_1; f_2; A_n^1) &\leq V_{\rm B}(g_*; f_2; A_n^1).
\end{align}
We combine (\ref{Vb2})--(\ref{Vb4}) to get
\begin{align*}
\mu(f_1) V_{\rm B}(\sigma f_1; f_2; A_n^1) &+(1-\mu(f_1)) V_{\rm B}(\phi f_1; f_2; A_n^1)\\
\leq & \mu(f_1) V_{\rm B}(\sigma \tilde{f}_1; f_2; A_n^1) + (1-\mu(f_1)) V_{\rm B}(\phi\tilde{f}_1; f_2; A_n^1). 
\end{align*}
Applying (\ref{Vb1}) then yields 
$$V_{\rm B}^1(f_1; f_2; A_n)\leq V_{\rm B}^1(\tilde{f}_1; f_2; A_n).$$ 
The rest of the proof is standard. 
\end{proof}

{\bf Remark.}  Theorem~\ref{thm2} focuses on the parameter $p$.  If we still require equal prior means for $p$, but impose the log-concavity order on $\theta=\log (p/(1-p))$ rather than $p$, then $V_{\rm B}$ is ordered by virtually the same proof.  This result is distinct from Theorem~\ref{thm2} because the relative log-concavity order is usually not preserved by monotone transformations. 

\section{Normal bandits with general priors}
The main result of this section (Theorem~\ref{thm3}) extends Theorem~\ref{thm1} to general priors for normal bandits.  Similar to Theorem \ref{thm2}, Theorem~\ref{thm3} is based on the relative log-concavity order, although it is more restrictive because we only compare a general prior with a normal prior. 

Given $\theta_i,\ i=1,2,$ let us assume that observations from arm $i$ are i.i.d.\ ${\rm N}(\theta_i, 1)$.  Priors on $\theta_i$ are independent with Lebesgue densities $f_i$.  We shall denote the value of this normal bandit with discount sequence $A_n$ by $V_{\rm N}(f_1; f_2; A_n)$.  Denote the mean of any $f$ by $\mu(f)=\int_{-\infty}^\infty \theta f(\theta)\, {\rm d}\theta$. 

\begin{theorem}
\label{thm3}
Let $\tilde{f}_1\equiv {\rm N}(\alpha, 1/\tau)$. 
\begin{enumerate}
\item
If $f_1\leq_{\rm lc} \tilde{f}_1$ and $\mu(f_1)= \alpha$, then $V_{\rm N}(f_1; f_2; A_n)\leq V_{\rm N}(\tilde{f}_1; f_2; A_n)$.
\item
If $\tilde{f}_1 \leq_{\rm lc} f_1$ and $\mu(f_1)= \alpha$, then $V_{\rm N}(\tilde{f}_1; f_2; A_n)\leq V_{\rm N}(f_1; f_2; A_n)$.
\end{enumerate}
\end{theorem}
Let $\Lambda_{\rm N}(f; A_n)$ denote the break-even value of a one-armed normal bandit with prior $f$ for the mean of the unknown arm.  We obtain Corollary~\ref{coro4} as a consequence of Theorem~\ref{thm3} and Lemma~\ref{lemlam}.
\begin{corollary}
\label{coro4}
Assume $A_n$ is regular and $a_1>0$.  Define $\tilde{f}\equiv {\rm N}(\alpha, 1/\tau)$. 
\begin{enumerate}
\item
If $f\leq_{\rm lc} \tilde{f}$ and $\mu(f)=\alpha$, then $\Lambda_{\rm N}(f; A_n)\leq \Lambda_{\rm N}(\tilde{f}; A_n)$.  
\item
If $\tilde{f}\leq_{\rm lc} f$ and $\mu(f)=\alpha$, then $\Lambda_{\rm N}(\tilde{f}; A_n)\leq \Lambda_{\rm N}(f; A_n)$.  
\end{enumerate}
\end{corollary}

The condition $f\leq_{\rm lc} {\rm N}(\alpha, 1/\tau)$ is essentially ${\rm d}^2\log f(\theta)/{\rm d}\theta^2 \leq -\tau$, which can be regarded as a strong form of information ordering.  The appearance of $\leq_{\rm lc}$ is therefore especially intuitive in Theorem~\ref{thm3} and Corollary~\ref{coro4}.  It is an open problem whether Theorem~\ref{thm3} and Corollary~\ref{coro4} hold without assuming that one of the priors is normal. 

The rest of this section proves Theorem~\ref{thm3}.  We need a technical result (Lemma~\ref{lem3}) which may be of independent interest. 

\begin{lemma}
\label{lem3}
Let $g$ be a differentiable function on $\mathbf{R}$.  Assume $X$ is a random variable satisfying $Eg(X)=EX$.  
\begin{enumerate}
\item
If $0\leq g'(x)\leq 1,\ x\in \mathbf{R}$, then $g(X)\leq_{\rm cx} X$.
\item
If $g'(x)\geq 1,\ x\in \mathbf{R}$, then $X\leq_{\rm cx} g(X)$.
\end{enumerate}
\end{lemma}
\begin{proof}
We prove Part 1 only.  Part 2 follows from Part 1 by considering the inverse function of $g$.  As $Eg(X)=EX$, one criterion for $g(X)\leq_{\rm cx} X$ is 
\begin{equation}
\label{cx}
E\max\{0, g(X)-b\}\leq E\max\{0, X-b\},\quad b\in \mathbf{R}.
\end{equation}
See, e.g., Shaked and Shanthikumar (2007; Theorem 3.A.1).  Let us assume $0\leq g'(x)\leq c$ for some $0<c<1$.  Otherwise we consider $cg(x)$ and let $c\uparrow 1$.  As $g(x)$ is a contraction, it has a unique fixed point, say $x_0$.  Consider two cases. 

Case (i): $b\geq x_0$.  If $x\geq x_0$ then $g(x)-g(x_0)\leq x-x_0$, i.e., $g(x)\leq x$, and $\max\{0, g(x)-b\}\leq \max\{0, x-b\}$.  If $x< x_0$ then $g(x)\leq g(x_0)=x_0$ and 
$$\max\{0, g(x)-b\}\leq \max\{0, x_0-b\}=0\leq \max\{0, x-b\}.$$
In either case $\max\{0, g(x)-b\}\leq \max\{0, x-b\}$, which implies (\ref{cx}). 

Case (ii): $b< x_0$.  Applying the argument of Case (i) to $\tilde{g}(x)\equiv -g(-x)$ and $\tilde{X}\equiv -X$ yields
$E\max\{0, b-g(X)\}\leq E\max\{0, b-X\}$, which reduces to (\ref{cx}) because $Eg(X)=EX$. 
\end{proof}

\begin{proof}[Proof of Theorem~\ref{thm3}]
We only prove Part 1; the second part is similar.  The $n=1$ case is easy.  For $n\geq 2$ we use induction.  The equations (\ref{Vdef1})--(\ref{Vdef3}) become 
\begin{align}
\nonumber
V_{\rm N}(f_1; f_2; A_n) &= \max\left\{V_{\rm N}^1(f_1; f_2; A_n),\ V_{\rm N}^2(f_1; f_2; A_n)\right\};\\
\label{VN1}
V_{\rm N}^1(f_1; f_2; A_n) & = a_1 \mu(f_1) + E\left[V_{\rm N}(f_1^X; f_2; A_n^1)|\Phi f_1\right];\\
\nonumber
V_{\rm N}^2(f_1; f_2; A_n) & = a_1 \mu(f_2) + E\left[V_{\rm N}(f_1; f_2^Y; A_n^1)|\Phi f_2\right]. 
\end{align}
We denote the posterior $f_1^x(\theta)\propto f_1(\theta) \exp[-(x-\theta)^2/2]$; similarly for $f_2^y$.  In $E[g(X)|\Phi f]$, the density of $X$, denoted by $\Phi f$, is the convolution of $f$ with the standard normal.  (Note the difference from the notation in Section~2.)  Let $m(x; f)$ denote the posterior mean of $\theta$ when $x$ is observed and the prior is $f$, i.e., $m(x; f)=\int_{-\infty}^\infty \theta f^x(\theta)\, {\rm d}\theta$.  Direct calculation yields 
\begin{equation}
\label{heat}
\frac{{\rm d} m(x; f)}{{\rm d}x} = Var(\theta|f^x).
\end{equation}
That is, the derivative of $m(x; f)$ is simply the posterior variance of $\theta$. 

Suppose $f_1\leq_{\rm lc} \tilde{f}_1\equiv {\rm N}(\alpha, 1/\tau)$ and $\mu(f_1)=\alpha$.  Then 
\begin{equation}
\label{normlc}
f_1^x\leq_{\rm lc} {\rm N}\left(m(x; f_1), \frac{1}{\tau+1}\right).
\end{equation}
It can be shown that (i) if $X$ is distributed as $\Phi f_1$, then $m(X; f_1)\leq_{\rm cx} (X+\tau \alpha)/(\tau +1)$; (ii) $\Phi f_1$ is smaller than $\Phi\tilde{f}_1\equiv {\rm N}(\alpha, 1+1/\tau)$ in the convex order.  To prove (i), note that (\ref{normlc}) holds with $\leq_{\rm lc}$ replaced by $\leq_{\rm cx}$ as the two sides have equal means.  By (\ref{heat}) we have 
$$0\leq \frac{{\rm d} m(x; f_1)}{{\rm d}x}\leq \frac{1}{\tau +1},\quad x\in \mathbf{R}.$$
If $X$ is distributed as $\Phi f_1$ then both $(X+\tau\alpha)/(\tau+1)$ and $m(X; f_1)$ have mean $\mu(f_1)=\alpha$.  Thus claim (i) holds by Lemma~\ref{lem3}.  Claim (ii) holds because $f_1\leq_{\rm cx} \tilde{f}_1$ and the convex order is closed under convolution. 

We have 
\begin{align*}
E\left[V_{\rm N}(f_1^X; f_2; A_n^1)|\Phi f_1\right] &\leq \left. E\left[V_{\rm N}\left({\rm N}\left(m(X; f_1), \frac{1}{\tau+1}\right); f_2; A_n^1\right)\right|\Phi f_1\right]\\
&\leq \left. E\left[V_{\rm N}\left(\tilde{f}_1^X; f_2; A_n^1\right)\right|\Phi f_1\right]\\
&\leq \left. E\left[V_{\rm N}\left(\tilde{f}_1^X; f_2; A_n^1\right)\right|\Phi\tilde{f}_1\right], 
\end{align*}
where the first inequality holds by (\ref{normlc}) and the induction hypothesis, the second by claim (i), noting
$$\tilde{f}_1^X={\rm N}\left(\frac{X+\tau\alpha}{\tau+1}, \frac{1}{\tau+1}\right),$$
and the third by claim (ii).  The last two inequalities also use the convexity of $V_{\rm N}$ with respect to the mean of a normal prior, i.e., Proposition~\ref{prop2}.  (Although Proposition~\ref{prop2} assumes normal priors for both arms, this can be relaxed.)  It follows from (\ref{VN1}) that 
$$V_{\rm N}^1(f_1; f_2; A_n)\leq V_{\rm N}^1(\tilde{f}_1; f_2; A_n).$$
The rest of the proof is standard.
\end{proof}


\section{Discussion}
Results in previous sections suggest the following conjecture.  Consider a two-armed bandit in the general exponential family setting with conjugate priors.  Suppose the prior expected yield of one pull from each arm is the same, but the prior weight of arm 1 is larger.  Then it seems reasonable that arm 2 is optimal at the first stage, i.e., in the notation of Section~3, 
$$\frac{\gamma_1}{\tau_1}=\frac{\gamma_2}{\tau_2}\quad {\rm and}\quad \tau_1>\tau_2\quad \Longrightarrow\quad \Delta(\gamma_1, \tau_1; \gamma_2, \tau_2; A_n)\leq 0.$$ 
This holds if the discount sequence is infinite-horizon geometric.  Indeed, it is optimal to pull arm 2 because, according to Corollary~\ref{coro1}, arm 2 has a larger Gittins index.  For non-geometric discounting, we cannot apply Corollary~\ref{coro1} due to the lack of an index policy.  In fact, Berry (1972) proposed this conjecture for Bernoulli bandits with uniform discounting, and this special case is still open. 


\end{document}